\documentclass[11pt]{article}
\usepackage{amsmath,amssymb,graphicx,float,epsf}
\def\RR{\hbox{I\kern-.2em\hbox{R}}}
\newcommand{\qed}{\hbox to 0pt{}\hfill$\rlap{$\sqcap$}\sqcup$ \vspace{3mm}}
\numberwithin{equation}{section}
\restylefloat{figure}
\usepackage{graphicx} 
\usepackage[table,xcdraw]{xcolor}
\usepackage{threeparttable}
\usepackage{amsmath}
\usepackage{amsthm}
\usepackage{amsfonts}
\usepackage{ mathrsfs }
\usepackage[caption = false]{subfig}
\usepackage{wrapfig}
\usepackage{enumerate}
\usepackage{rotating}
\usepackage{enumitem}
\usepackage{array}
\usepackage{longtable}
\usepackage{subfig}
\usepackage[margin=1cm]{caption}
\usepackage[margin=1in]{geometry}
\usepackage{authblk}

\usepackage[utf8]{inputenc}

\newtheorem{Th}{Theorem}
\newtheorem{Cor}{Corollary}

\newtheorem{Conj}{Conjecture}
\usepackage{graphicx}
\usepackage{hyperref}
\hypersetup{
    colorlinks=true,
    linkcolor=blue,
    filecolor=magenta,
    urlcolor=cyan,
    citecolor=blue,
}
\date{}
\graphicspath{ {Images/} }
\begin{document}

    \title{Cryptanalysis of RSA Cryptosystem: Prime Factorization using Genetic Algorithm}

    \author[1,2]{\small Mahadee Al Mobin\thanks{Email: mahadeealmobin@gmail.com}}
    \author[2]{\small Md Kamrujjaman\thanks{Corresponding author email: kamrujjaman@du.ac.bd}}

    \affil[1]{\footnotesize Bangladesh Institute of Governance and Management, Dhaka 1207, Bangladesh}
    \affil[2]{\footnotesize Department of Mathematics, University of Dhaka, Dhaka 1000, Bangladesh}

    \maketitle

    \vspace{-1.0cm}
    \noindent\rule{6.35in}{0.02in}\\
    {\bf Abstract.}
    Prime factorization has been a buzzing topic in the field of number theory since time unknown. However, in recent years, alternative avenues to tackle this problem are being explored by researchers because of its direct application in the arena of cryptography. One of such applications is the cryptanalysis of Rivest–Shamir–Adleman (RSA) numbers, which requires prime factorization of large semiprimes. Based on numerical experiments, this paper proposes a conjecture on the distribution of digits on prime of infinite length. This paper infuses the theoretical understanding of primes to optimize the search space of prime factors by shrinking it upto $ 98.15\% $, which, in terms of application, has shown $26.50\%$ increase in the success rate and $41.91\%$ decrease of the maximum number of generations required by the genetic algorithm used traditionally in the literature. This paper also introduces a variation of the genetic algorithm named ``Sieve Method" that is fine-tuned for factorization of big semi-primes, which was able to factor numbers up to 23 decimal digits with $84\%$ success rate. Our findings shows that sieve methods on average has achieved $321.89\%$ increase in success rate and $64.06\%$ decrement in the maximum number of generations required for the algorithm to converge compared to the existing literatures.\\

    \noindent
    \noindent\rule{6.35in}{0.02in}

    \section{Introduction}\label{sec:introduction}
    Rivest-Shamir-Adleman (RSA) numbers are crucial for the security of RSA encryption, a foundational technology for secure digital communication.
    The main goal of cryptography is to keep adversaries like Oscar from deciphering the information being sent while ensuring secure communication between two parties, usually Alice and Bob, across an unsecure channel. The methodical process of decrypting ciphertext into plaintext without having access to the decryption key is known as cryptanalysis. This endeavor is in line with the goals of an adversary, represented by the figure of Oscar, who is attempting to intercept and interpret the messages that Alice and Bob are exchanging. Symmetric key cryptography relies on prior agreement and stringent confidentiality between the sender and the receiver because they use a single shared key for both encryption and decryption. On the other hand, asymmetric key cryptography use distinct keys for encryption and decryption, where the public key is used for encryption and the private key is used for decryption. Asymmetric key systems, mostly rely on computational security, in which the decryption algorithm is purposefully difficult but the encryption algorithm is readily computed, a construction analogous to that of a trapdoor. The message's confidentiality is underscored by the fact that only those with the private key can decrypt it.\\

    One of the principal issues in the field of asymmetric key cryptography is integer factorization, which has fascinated scholars for millennia. This cryptographic system protects the privacy of messages by making sure that only the intended receiver is able to decode them. A prominent example of this protocol is the RSA (Rivest-Shamir-Adleman) cryptosystem \cite{rivest1978method}. The security that supports RSA arises from the fact that multiplying large prime numbers is inherently simple, while factoring huge semi-primes numbers is a daunting challenge \cite{rivest1978method} as its NP-Completeness has not been established.\\

    Being an asymmetric cryptographic scheme, RSA has two keys namely:
    \begin{enumerate}[]
    	\item Public Key: $ (M,s) $
    	\item Private Key: $ (j,k,r) $
    \end{enumerate}

    where, $ M $ is a huge semi-prime constructed using the primes $ j $ and $ k $. Thus, Euler's totient function, $ \phi(M)=(j-1)(k-1) $. $ r $ and $ s $ are integers such that $ rs\equiv 1 (mod\phi(M)) $ where $ 1<s<\phi(M) $.\\

    If $ p $ is the plain text and $ c $ is the ciphered text then, we can present the rule to encrypt and decrypt in the following fashion:
    \begin{eqnarray*}
    	e_k(p)=p^s (mod M)\\
    	d_k(c)=c^r (mod M)
    \end{eqnarray*}

    The formulae above point out an important missing element from the decryption rule i.e., the number $ r $. There are multiple methods for figuring out this missing piece when it comes to RSA encryption. Calculating Euler's totient function $ \phi(M) $ is one approach. On the other hand, factorizing the semi-prime $ M $ to find its prime factors $ j $ and $ k $ is the most common attack against RSA. Euler's totient function may then be calculated from these integers, making it easier to determine $ r $, considering the fact that $ s $ is publicly available.\\

    The pursuit of solving the integer factorization problem has produced a multitude of different approaches in recent years. According to previous research \cite{mishra2016heuristic,yampolskiy2010application}, these methodologies can be categorized into three principal classes: specialized, generic, and alternative. Specialized strategies leverage specific numerical structures to optimize their effectiveness. On the other hand, the generic method focuses on developing integer factorization solutions that can be used to all kinds of numbers. On the contrary, the alternative strategy departs from traditional approaches and looks for new ways to handle the issue. This frequently involves using computational intelligence methods, such as particle swarm intelligence, firefly algorithms, neural networks, and genetic algorithms (GA).\\

    The field of cryptanalysis has experienced the use of a wide range of computational intelligence techniques, most of which have been used to classical ciphers—those that are considered insecure and therefore not often used in modern settings. Famous examples include traditional ciphers like the substitution cipher \cite{mudgal2017application}, the Purple cipher \cite{shikhare2015cryptanalysis}, Substitution Permutation Networks (SPN) \cite{brown2009genetic}, the RC4 Stream cipher \cite{ferriman2014solving}, and the Tiny Encryption Algorithm (TEA) \cite{ma2011evolutionary}. By significantly increasing the complexity against theoretical attacks, the application of GAs to the RC4 encryption shows significant success \cite{ferriman2014solving}. Studies pertaining to the application of GA on Tiny Encryption Algorithm have also provided evidence indicating that keys with a higher proportion of random words or that are composed only of one-bit words are more resilient to possible assaults \cite{ma2011evolutionary}.\\

    In \cite{brown2009genetic}, researchers used a genetic algorithm to find weak keys in the Substitution Permutation Network cipher. Their findings suggested that the use of evolutionary algorithms in cryptanalysis could be critical for detecting weak keys and therefore improving cipher security. This argument is based on the understanding that anyone with decent computational resources can use genetic algorithms to potentially undermine ciphers if weak keys are used. Furthermore, the paper urged for a greater application of genetic algorithms to non-classical ciphers, arguing that such efforts could shed light on their efficacy in strengthening cryptographic systems \cite{brown2009genetic}.\\

    Similarly, recent studies have used genetic programming (GP) to improve the cryptanalysis of elliptic curve cryptosystems \cite{ribaric2017genetic}. These cryptosystems, which comprise public-key infrastructure, rely on the complexity of the elliptic curve discrete logarithm issue. Within the scope of the study, genetic programming aided in the speeding of a component of a well-established algorithm designed to solve this difficult challenge.\\

  	Recently, researchers have begun to investigate bio-inspired algorithms as potential solutions to integer factorization. Notably, studies such as those cited in \cite{chan2002automatic}, \cite{chang2005fast}, and \cite{meletiou2002first} investigated the use of alternative approaches, such as genetic programming (GP) and neural networks, to address the difficulty of integer factorization.\\

    \cite{yampolskiy2010application} used genetic algorithm to solve the problem of integer factorization. The algorithm used a chromosomal representation to encapsulate the two prime components, indicated as $ p $ and $ q $, which generate the integer $ N $ via the equation $ N = p \times q $. The chromosomal length was mapped to the decimal representation of $ N $. The chromosome was partitioned in accordance with the premise that each prime factor $ p $ and $ q $ was no longer than half the length of $ N $. The first half of the chromosome represented the decimal value of $ p $, while the second half depicted the decimal value of $ q $, as shown in equation \eqref{eq:ga_chromosome_presentation}.

    \begin{equation}\label{eq:ga_chromosome_presentation}
    	[p_1p_2p_3\dots p_{|N|/2}q_1q_2q_3\dots q_{|N|/2}]
    \end{equation}

    The fitness function was designed to assess the similarity between the product generated by multiplying the prime components $ p $ and $ q $ taken from the chromosome and the number $ N $, using parity as the metric. In the case where $ N $ had $ m $ decimal digits, the optimal fitness value assigned to the chromosome equaled $ m $. This situation represented an exact match between all $ m $ digits of the chromosome-derived product and $ N $, indicating that the proper $ p $ and $ q $ had been identified. However, the occurrence of several local minima presented a hurdle, especially when the product formed from $ p $ and $ q $ differed by a single digit from $ N $. Notably, the study factorized a twelve-digit semi-prime, 103694293567 = 143509 * 722563, which required little more than six hours.\\

    Mishra et. al.  uses multi-threaded bound changing chaotic firefly technique, for the integer factorization problem \cite{mishra2014multithreaded}. The firefly algorithm is modeled after the collective behavior of fireflies. It uses a fitness function to lure fireflies closer each other according to their light levels. Ten different test datasets were used to assess the algorithm's performance. The biggest integer that was examined was 51790308404911 (5581897 * 9278263), which is 14 digits (or 46 bits) in length. Interestingly, the method factored this figure with an accuracy of between 80\% and 100\%, depending on how many fireflies were used in the computation.\\

    In order to tackle the problem of integer factorization, \cite{mishra2016heuristic} utilized a heuristic approach based on molecular principles. Motivated by the atomic structure in a space with negligible interatomic forces, the algorithm combines an energy function, a force function, and a movement function to determine possible solutions. Using the same biggest integer tried in \cite{mishra2014multithreaded}, 51790308404911 (5581897 * 9278263), the algorithm's effectiveness was evaluated. It is noteworthy that the algorithm factored this number with a remarkable success rate of 69\%. Additionally, the authors performed a comparative analysis by comparing the results of their system to a random search algorithm. The 11-digit (35-bit) number 42336478013 was notably factorable by the random search method, but with a meager success rate of 4\%.\\

    Three genetic algorithms were used in \cite{rutkowski2020cryptanalysis}; the most successful approach factorized values up to 22 decimal digits with a success rate of $ 3.33\% $ by using a chromosome to represent $ m $ in the equation $ prime = 6m \pm 1 $.The best algorithm was able factorize 15 digits number ($ 115137038087959=10037141\times 11471099 $) with a success rate of $ 36.67\% $. This accomplishment represented a significant improvement over earlier approaches. However, the inability of the proposed genetic algorithms to factor integers with key lengths longer than fifteen is evident. Moreover, RSA factorization's difficulty with higher integers presents a possible obstacle that can lower the algorithm's success rate.\\

    Thus the computational intelligence approaches observed from the review of literatures, shows that the best achievement presented by these class of algorithms in integer factorization problems in a stable fashion was 15 digits number ($ 115137038087959=10037141\times 11471099 $) with a success rate of $ 36.67\% $ by a GA namely ``Chromosome is m''. The algorithm was also able to factorize a 22 decimal digits luckily . \emph{This paper aims to improve over the current body of literature to push the capabilities of the GA using theoretical constructs of RSA and prime numbers to be able to factorize bigger RSA numbers in a much stable fashion with higher success rate. }  \\

    This paper adopts alternative methodological approach which attempts to factorize $ M $  using GA jeopardizing  the public key cryptosystem's security. The paper exhibits two variations of GA namely: Simple GA which is modified adaptation of the Simple GA works presented in \cite{yampolskiy2010application,rutkowski2020cryptanalysis} with better search space optimization and ``Sieve Method'' which is theoretical generalization of the ``Chromosome is m'' GA presented in \cite{rutkowski2020cryptanalysis}. The two GA variants that are explained in this paper are significantly better than the current algorithms that are reported in the literature in terms of computational efficiency, ability to factorize the bigger numbers, and success rates. \\

    The value addition of the paper to the current body of literature is illustrated as follows:
    \begin{enumerate}
        \item The paper aims to propose a conjecture on the distribution of digits on larger primes.

        \item The aim to propose to GA variants that significantly better than the current algorithms that are reported in the literature in terms of computational efficiency, ability to factorize the bigger numbers, and success rates.

        \item The paper presents improved strategies for initializing the search space, based on theoretical insights into RSA and prime numbers. These changes resulted in a significant reduction in the search space of up to $ 98.15\% $, adding to the increased robustness of these bio-inspired approaches.

    \end{enumerate}

    The remaining segment of the paper is structured as follows. Section \ref{sec:theorybackground} explores the existing and novel theoretical  development regarding cryptography, RSA, semi-primes and primes. Section \ref{sec:method} discusses the simple genetic algorithm and a variation of genetic algorithm named ``Sieve Algorithm" that is proposed in this article. Section \ref{sec: specs} presents the technical specifications of the software and hardware used for the cryptanalysis. Section \ref{sec:result} presents our finding using each of the discussed algorithms and does comparative analysis with the existing literatures based on key performance indicators (KPIs). Section \ref{sec:lim} discusses the further avenues that can be explored in the scope of this study. Section \ref{sec:conclusion} finally presents the conclusion.

    \section{Theoretical Background}\label{sec:theorybackground}

	We use some theoretical developments on the RSA and prime numbers to optimize the search space and candidate selection criterion for our GAs. There are mainly two theorems that we use to make our algorithm more robust.\\

	The first theorem that we use is a work by \cite{xingbo2019number} which through a series of theoretical development proves the following result:

	\begin{Th}\label{theorem:digits_of_factor}
		If M is RSA number such that, $ M=j\times k $, $ 1<j<k $ and $ D_j $ and $ D_k $ represent the length of the prime factors then,
		\begin{equation}\label{eq:divisorlength}
			 D_j=D_k=\left\{\displaystyle\begin{array}{cc}
			\frac{D_M}{2} &, D_M \text { is even } \\
			\frac{D_M+1}{2}&, D_M \text { is odd }
		\end{array}\right.
		\end{equation}
	\end{Th}

	This theorem gives us an idea about the length of the prime divisors of a RSA number hence has been a game changer for search space optimization. Previously, the search space considered for factors of M is $ (2, \sqrt{M}) $ but now based on this theorem the search space for the factor of M can considered as $ (10^{D_j-1},\sqrt{M}) $.
	 Now, the percentage of shrinkage in solution space can be expressed as \begin{equation}\label{eq:ss_shrinkage}
		 \text{Search Space Shrinkage (SSS)}=\displaystyle \frac{10^{D_j-1}-2}{\sqrt{M}-2}\times 100
	\end{equation}\\

	We have calculated the shrinkage of solution space for different length of numbers to have an idea of the actual optimization achieved due to this theoretical development. Our findings are illustrated in figure \ref{fig:ss_optim_line}.\\

	\begin{flushleft}
		\includegraphics[width=\textwidth]{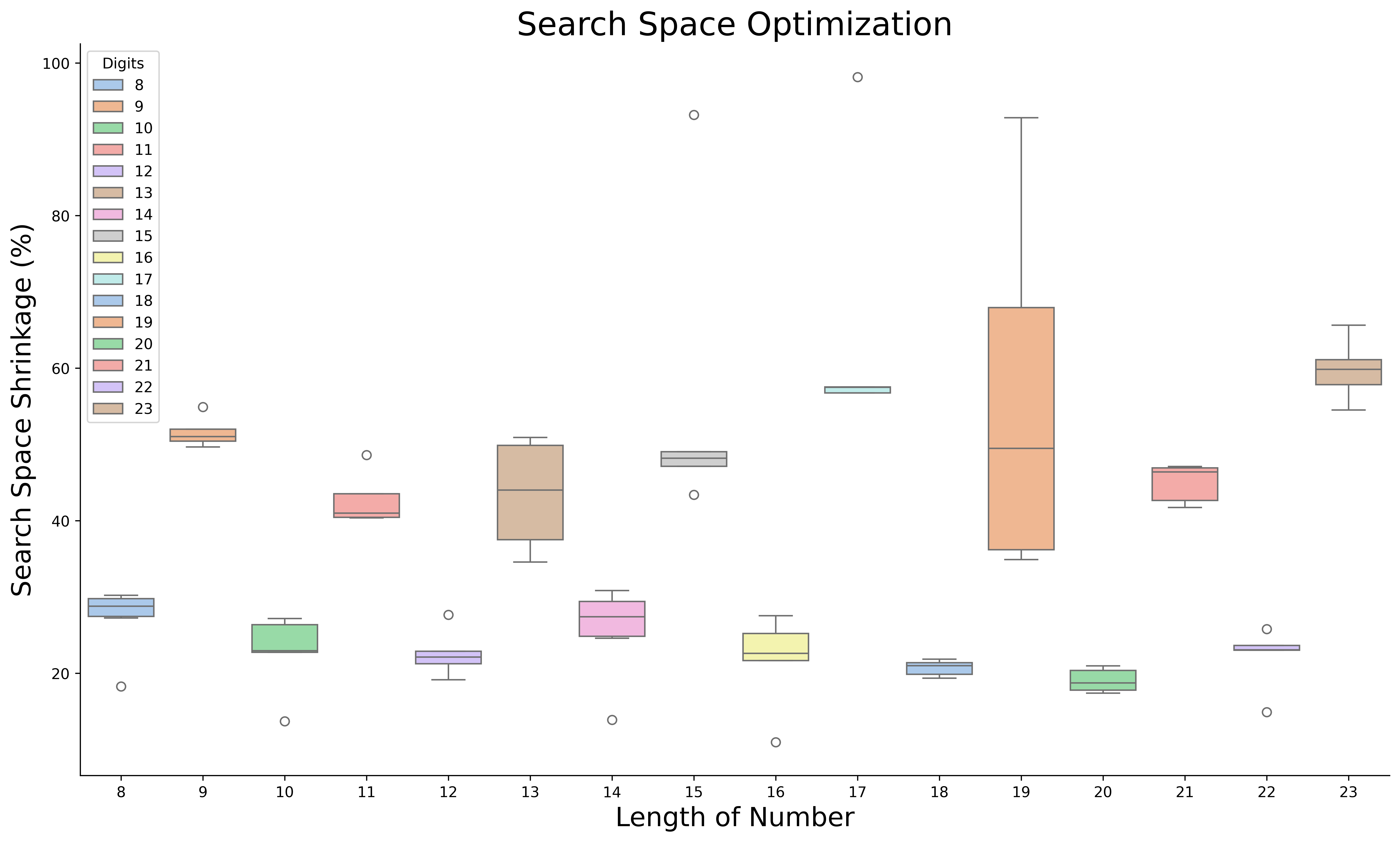}
		\captionof{figure}{Search space optimization of GA using theoretical development for varying length of numbers.  }\label{fig:ss_optim_line}
		\footnotesize{The figure depicts the shrinkage of the search space in percentage for varying length of RSA numbers. The simulation has been conducted on the five datasets that we considered for this paper each of which contains a RSA number of length between 8 and 23. Due to this theoretical development, we were able to shrink the search space upto $ 98.15\% $ based on the simulation done on the 81 unique RSA numbers.}
	\end{flushleft}
	The figure depicts the variation of SSS obtained by the theoretical development for the five datasets of RSA numbers we have considered for cryptanalysis in this paper. Details on the datasets can be found in appendix \ref{app:sieve_algo_output} and \ref{app:simple_ga_output}. The figure shows that we were able to achieve up to $ 98.15\% $ SSS. The graph shows an oscillatory patter for SSS which suggest that SSS shows variation for odd and even length of numbers. To investigate the distribution of SSS and the oscillatory pattern, we draw an overall box plot of the SSS and a grouped boxplot for odd and even length of numbers which are illustrated in figure \ref{fig:ss_optim_box} and \ref{fig:ss_optim_box_grouped} respectively.\\

	\begin{flushleft}
		\includegraphics[width=\textwidth]{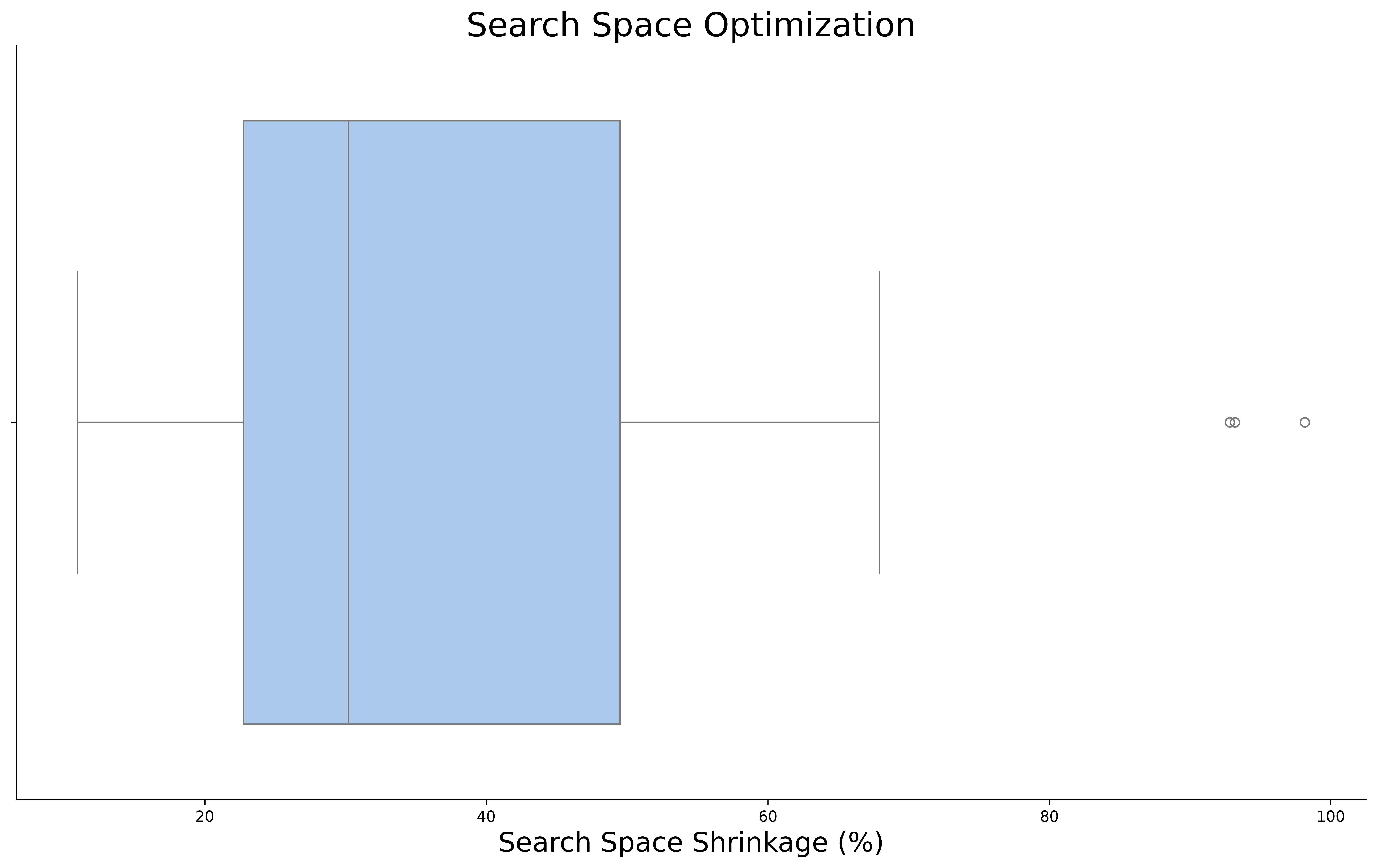}
		\captionof{figure}{Distribution of SSS obtained in the experimentation. }\label{fig:ss_optim_box}
		\footnotesize{The figure depicts the distribution of SSS. Due to this theoretical development, we were able to shrink the search space upto $ 98.15\% $ with a median value of $ 30.22\% $.}
	\end{flushleft}

	\begin{flushleft}
		\includegraphics[width=\textwidth]{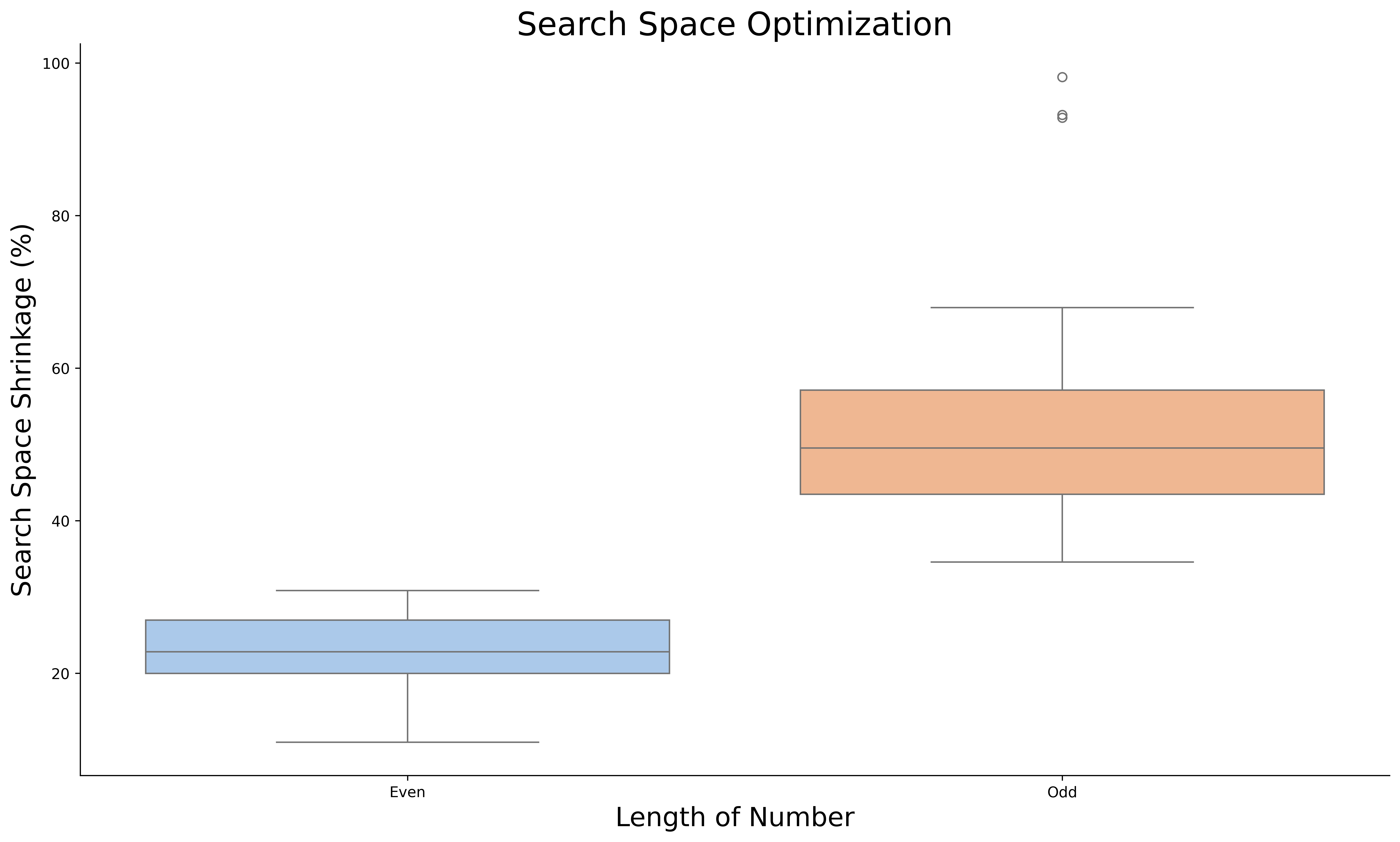}
		\captionof{figure}{Distribution of SSS obtained in the experimentation grouped by the length of the number. }\label{fig:ss_optim_box_grouped}
		\footnotesize{The figure depicts the distribution of SSS for odd and even length of numbers. In case of numbers with even length, we achieved SSS upto $ 30.86\% $ with a median value of $ 22.82\% $ and for that of number with odd length, we achieved SSS upto $ 98.15\% $with a median value of $ 49.54\% $.}
	\end{flushleft}

	Figure \ref{fig:ss_optim_box} shows the overall distribution of SSS obtained in our experimentation. We were able to obtain up to $ 98.15\% $ SSS with a median value of $ 30.22\% $. To further investigate the skewness of the distribution, we drew the box plot of the distribution of SSS grouped by the nature of the length of the number as illustrated in figure \ref{fig:ss_optim_box_grouped}. It shows huge deviation observed between the two groups. We were able to higher SSS for the odd group than that of the even group. For the odd group, we achieved  SSS up to $ 98.15\% $ with a median value of $ 49.54\% $ whereas for their even counterpart, we achieved SSS up to $ 30.86\% $ with a median value of $ 22.82\% $. The reason for such deviation is evident from the equation \eqref{eq:divisorlength} and \eqref{eq:ss_shrinkage}. Equation \eqref{eq:divisorlength} which states that a RSA number of an even length, $ n $ and a RSA number of odd length, $ n-1 $ has factors of same length. Equation \eqref{eq:ss_shrinkage} states, that and  SSS is directly proportional to $ 10^{D_j}-1 $ and is inversely proportional to $ \sqrt{M}-2 $. This construction leads to bigger search space for even RSA numbers in comparison to it's preceding and following odd counterpart.\\

	In conclusion, using theorem \ref{theorem:digits_of_factor}, we were able to shrink the search space upto $ 98.15\% $ with a median value of $ 30.22\% $ and SSS is higher for odd RSA numbers then that of the even ones. This is major development over the existing literature in terms of computational time, memory and efficiency especially for bigger RSA numbers which is well illustrated in the result section of this paper.\\

	The final theorem that we use in the Dirichlet prime number theorem which states that,

	\begin{Th}[Dirichlet Prime Number Theorem]\label{theorem:dirichlet}
		$ \forall n\in \mathbb{N} $, there are infinitely many primes of the form $ an\pm d $ where $ (a,d)=1 $.
	\end{Th}

	A well known corollary of theorem \ref{theorem:dirichlet} is the following  upon which the ``Chromosome is m" algorithm is based.

	\begin{Cor}\label{cor:6m+-1}
		Every prime greater than 3 can be expressed in the form $ 6n\pm1 $.
	\end{Cor}

	We will use Theorem \ref{theorem:dirichlet} to construct the generalization of ``Chromosome is m" i.e. the ``Sieve Method".\\

	Finally, based on our experimentation on the probability distribution of each digit in the construct of a big prime numbers to device efficient mutation operation for decimal chromosome representation in GA, we propose a novel conjecture on the distribution of digits on prime which is stated below:
	\begin{Conj}[\emph{Mahadee-Kamrujjaman Conjecture}]
		Digits of a prime are uniformly distributed as the length of primes tends to infinity.
		 \begin{center}
		 	i.e.
		 \end{center}
		Let $ P_n (d) $be the probability of occurrence of a digit, $ d\in\{0,1,2,\dots,9\} $  in primes of length less
		than $ 10^n $  then for any choice of d,
		$$ \displaystyle \lim_{n\rightarrow \infty} P_n(d)=0.1$$
	\end{Conj}

	The experimentation result for the conjecture is illustrated in the figure \ref{fig:digitdistr}.

	\begin{flushleft}
		\includegraphics[width=\textwidth]{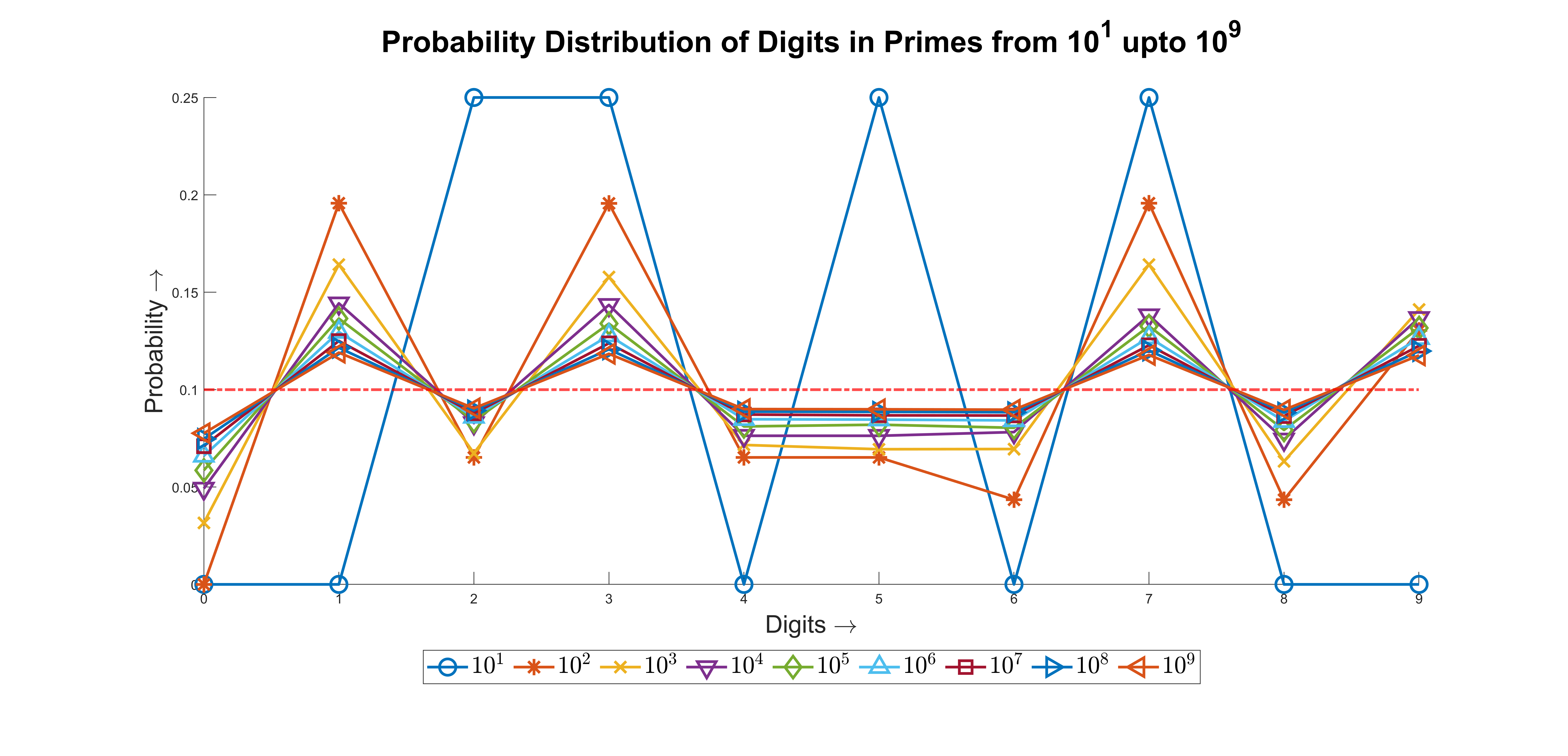}
		\captionof{figure}{Probability Distribution of Digits in Primes upto $ 10^9 $. }\label{fig:digitdistr}
		\footnotesize{The figure depicts the probability distribution of digits in a prime less than $ 10^n $ as $ n\rightarrow \infty $. Each of the curve in the figure depicts the probability distribution of digits in a prime less than $ 10^n $ a certain value of n. As the value of n increases (i.e. as the primes get bigger in length), we observe that the probability of occurrence of any digits moves closer to 0.1 (The red dotted line) irrespective of their probability for smaller values of n.}
	\end{flushleft}
	To provide proper justification for the conjecture, we device the following analytical approach for each iteration of n:
	\begin{enumerate}
		\item  We determined the primes less $ 10^n $.
		\item  We determine the total number of times each digit occurs in all the primes hence determine the probability of their occurrence.
		\item We plot the curve  of the probability distribution of the digits.
	\end{enumerate}
	Following the aforementioned approach, we were able to generate Figure \ref{fig:digitdistr}. Each of the curve in the figure represent the probability distribution of digits in a prime less than $ 10^n $ a certain value of n. As the value of n increases (i.e. as the primes get bigger in length), we observe that the probability of occurrence of any digits moves closer to 0.1 (The red dotted line) i.e. the probability of occurrence of all the digits become more or less uniform and it seems the behaviour will hold as n tends to infinity hence the formulation of the conjecture. Due to the limitation of the computational capacity, we were only able to validate our conjecture up to one billion. While our experimentation revealed the inefficiency of the decimal representation of the chromosome compared to its binary counterpart in terms of success rate and speed which has also been stated in \cite{rutkowski2020cryptanalysis}, its exploration nonetheless facilitated the discovery of a remarkable conjecture pertaining to prime numbers.\\
    \section{Method}\label{sec:method}

    This section discuss the two proposed algorithms namely Simple Genetic Algorithm (GA) and Sieve Method. Although the algorithms are different by construction, but both of these have some hyper parameters in common which are kept uniform throughout the algorithms. We consider the population 1500 for each generations and the algorithms runs upto 2000 generations until an optimized solution is achieved. Both the algortihms consider tournament selection size of 5.

    \subsection{Simple Genetic Algorithm}\label{subsec:simpleGA}

    Simple GA is the modification over the work  \cite{yampolskiy2010application,rutkowski2020cryptanalysis}. We have brought about significant changes in the construction which will be discussed in the following discussion.

    \subsubsection*{Initial Population}
   Rutkowski and Houghten   consider the initial population to be randomly generated 2000 between the interval $ (2,\sqrt{M}) $ \cite{rutkowski2020cryptanalysis}. We have modified the consideration search space to $ (10^{D_j-1},\sqrt{M}) $ as discussed previously. Following this development, we consider 1500 unique candidates randomly picked from the aforementioned search space.

    \subsubsection*{Chromosome Representation}
    Yampolskiy  used decimal representation to represent candidates \cite{yampolskiy2010application}. Further improvements were achieved by \cite{rutkowski2020cryptanalysis} using binary representation of candidates which considered the length of binary representation of candidate to be exactly half the binary representation of the $ M $. Rutkowski and Houghten  considered the first bit to be 1 to tackle the optimizer from converging to zero. Since a candidate must be from the interval $ (10^{D_j-1},\sqrt{M}) $, hence we consider the length of $ \lfloor \sqrt{M} \rfloor$ to be the length of each candidate for a problem \cite{rutkowski2020cryptanalysis}. We impose the condition for a candidate to be on the interval $ (10^{D_j-1},\sqrt{M}) $ as prerequisite to be included in a generation which frees our construction from the zero convergence problem.

    \subsubsection*{Fitness Function}
    In GA, fitness function is a function that summarizes how ``good" of a solution a candidate is for the problem at hand. In this context, the optimal solution is defined as the candidate solution that yields a remainder of zero, and a lower remainder signifies a superior solution. Therefore, we designate equation \eqref{eq:fitnessfunc} as the fitness function, representing the reciprocal of the remainder function, albeit with a removable discontinuity at zero.
    
    \begin{equation}\label{eq:fitnessfunc}
    	\text{f(p)}=\left\{ 
    		\begin{array}{cc}
    			\displaystyle\frac{1}{M\,modp}&\text{, M modp $ \neq 0$ }\\
    			K&\text{, M modp $ = 0$}
    		\end{array}\right.
    \end{equation} 
	Here, p is candidate solution for the problem, M is the large semi-prime and K is an absurdly large number. For our experimentation, we considered $ k=9.9\times 10^{37} $. As the fitness function by construction is reciprocal to that of the remainder, hence the for interpretation of fitness function in this regard we can state ``The higher, the better".\\
	
	To get an idea of how the fitness of the candidates are distributed in the search space with respect to the new fitness function, we plotted the fitness of each candidate in a log scale vs. the search space of the number $ 10909343 $ as illustrated in figure \ref{fig:solspace}.
    \begin{flushleft}
    	\includegraphics[width=\textwidth]{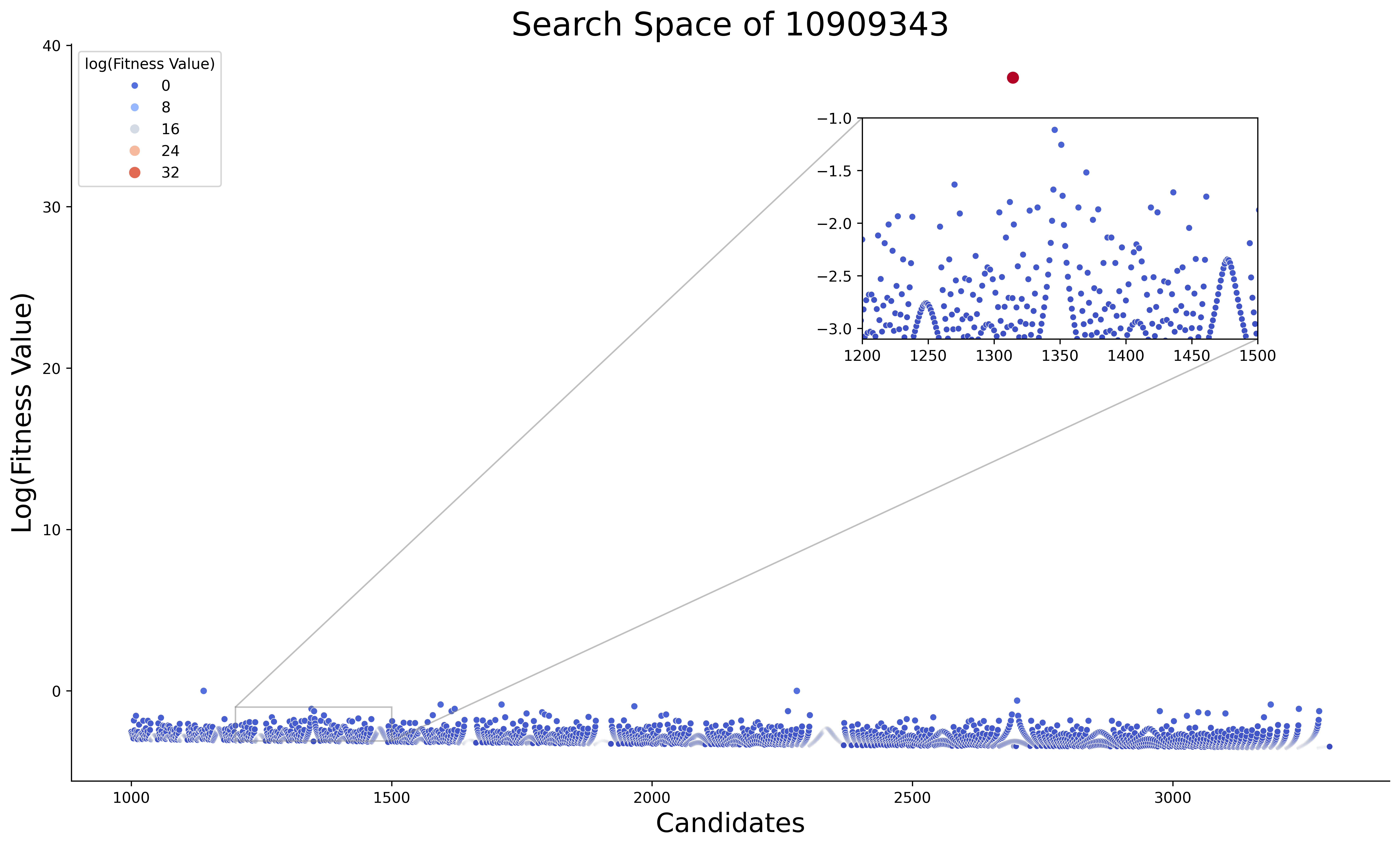}
    	\captionof{figure}{Solution Space of $ 10909343=2693\times4051 $ in Simple Genetic Algorithm. }\label{fig:solspace}
    	\footnotesize{The figure depicts the solution space of the number 10909343. We can easily detect the solution by the distinctively large value of the fitness function from the plot. The fitness function does not give negative value (By construction). As the y-axis is in log scale and almost all the fitness values are below 1 hence most of the values lie in the range of $ [-3,-1] $. From the zoomed section of the graph, it is evident that the solution space has multiple local maxima hence the algorithms in the paper shall adopt higher than usual mutation rate to combat the issue at hand. }
    \end{flushleft}
	
	The red dot in the figure depicts the solution to the problem. Upon closer inspection of the solution space, we observe that it is populated with a plethora of local maximas which makes it hard for GA to  converge to the actual solution. This problem is mitigated by adopting higher mutation rate
	\subsubsection*{Crossover}
	The paper adopts uniform crossover method to generate offsprings for future which randomly selects bits from the chromosomal representation of the parents to generate two offsprings.

	\subsubsection*{Mutation}
	The GA incorporates mutation through a stochastic process wherein a random locus within the chromosome is selected for modification. Upon selection, the value at the chosen locus undergoes alteration: a binary value of 0 transitions to 1, while a value of 1 transforms into 0.

    \subsection{Sieve Method}\label{subsec: sieve methods}
    
    The proposed method is a generalization of the ``Chromosome is m" method proposed in \cite{rutkowski2020cryptanalysis}. The proposed method is based on Theorem \ref{theorem:dirichlet} which considers a candidate solution, $ j $ of the form
    \begin{equation}\label{eq:sieve_method}
    	 j=an\pm d 
    \end{equation} 
    where, $ (a,d)=1 $. ``Chromosome is m" operates under the strict consideration of $ a=6 $, and $ d=1 $ as delineated in corollary \ref{cor:6m+-1}. Hence, the proposed method is a generalization over the aforementioned method. In this method, we try to find the optimal value of $ n $ using GA for equation \eqref{eq:sieve_method} such that the candidate $ j $ factorizes the large semi-prime $ M $. In case of simple GA, our search space for $j$ would be $ (10^{D_j-1},\sqrt{M}) $ but in the case of sieve algorithm, we have to optimize for the value of $ n $ using the equation \eqref{eq:sieve_method} hence our search space would be $ (\displaystyle \frac{10^{D_j-1}\mp d}{a},\displaystyle \frac{\sqrt{M}\mp d}{a}) $ which has been shrunk by a factor of $a$. The higher the value of $a$, the smaller the length of the search space i.e. the length of the search space is inversely proportional to the magnitude of $a$. This construction in unison with the search space optimization that has been demonstrated in this paper has shown better result in terms of performance which has discussed as length in the result section. To hyper parameter tune the value of $a$ and $d$, we consider $a\in [6,1000]\cap \mathbb{Z}$ and $d\in\{1,2\}$ and execute a GA (3 generations each of 25 candidates) over a small scale of the sieve method (200 generations each of 3000 candidates) only for numbers of length greater then or equal to nineteen. 
 
    \subsection*{Obstruction Based Learning}
     To increase variety in the candidate pool and hence improve algorithmic outcomes, the strategy includes computing the complement of each candidate relative to the search space for each generation, as shown by the following formula:
    
    \begin{equation}\label{eq:obl}
    	x^c=u+l-x=\sqrt{M}+10^{D_j-1}-x
    \end{equation}
    where, $ x $, $ x^c $ denotes a candidate and it's complement with respect to the search space respectively. $ u,\,l $ denotes the upper and lower bound of the search space respectively. This almost doubles the number of candidates considered for each generation and helps to address the scenario where the required solution lies on the peripheral of the search space. We strictly maintain that the uniqueness of a candidate for each generation by removing redundancies among the candidates.
    
    \subsubsection*{Fitness Function}
    The fitness function remains consistent with that of the Simple GA, as delineated in equation \eqref{eq:fitnessfunc}. However, additional preprocessing steps are undertaken prior to the evaluation of fitness.
    
    \subsubsection*{Crossover and Mutation}
    The crossover and mutation scheme adopted in this algorithm is same as that of the simple GA.

    \section{Technical Specification}\label{sec: specs}
	Both the algorithms are tested against five different datasets of RSA numbers which in total contains 81 RSA numbers. For each number, we execute 30 instance of  sieve method hence there are total $81\times30=2430$ instances of sieve methods to be executed. Similarly, there are 1530 instances of simple GA that are to be executed. The total computational time required for this instances is $1003872.94\text{ seconds}\equiv\text{11 days, 14 hours, 51 minutes, and 12.94 seconds}$.\\
	
	     \begin{flushleft}
	 	\includegraphics[width=\textwidth]{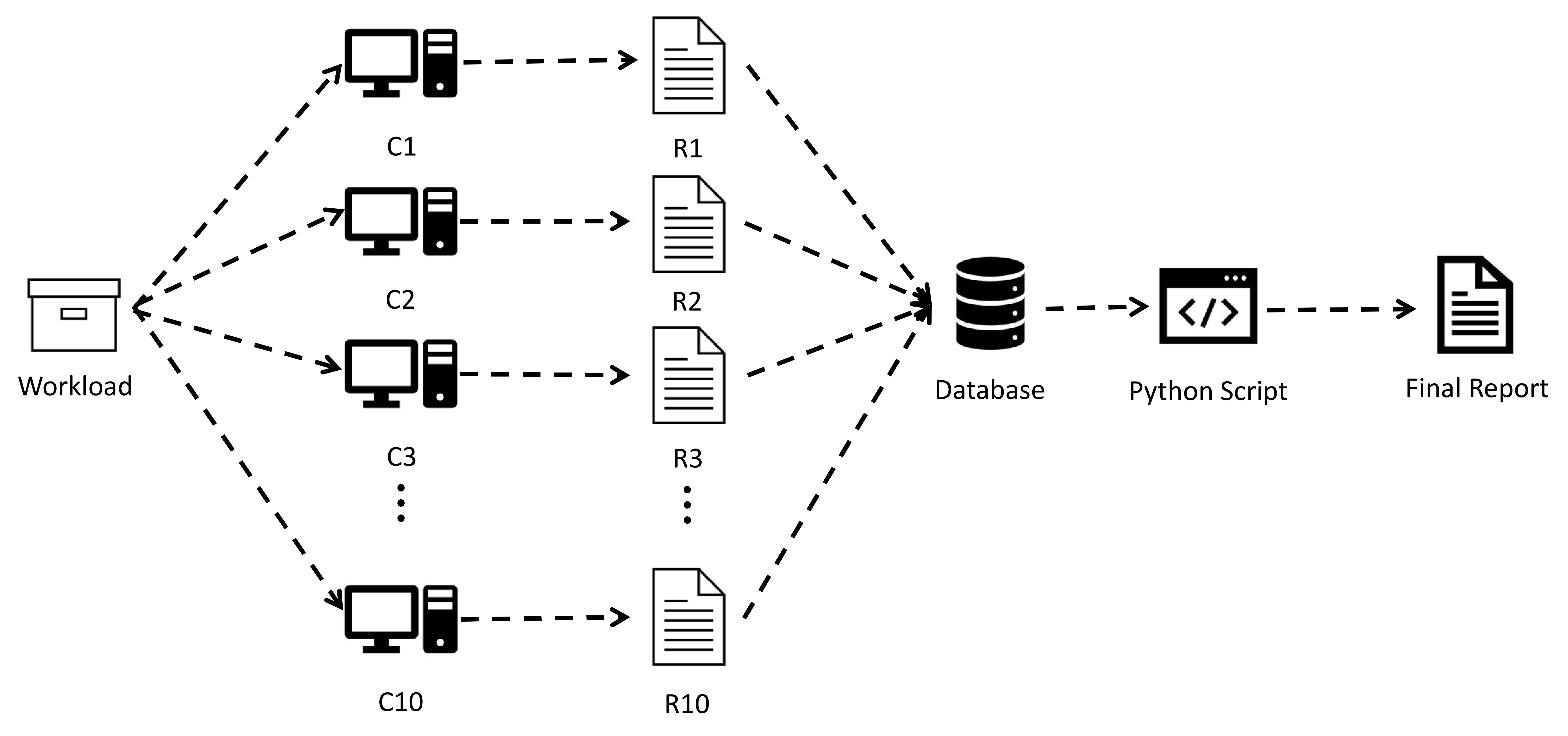}
	 	\captionof{figure}{Workflow of the simulation.}\label{fig:workflow}
	 	\footnotesize{The figure depicts the workflow of the entire simulation. Here, we used a fleet of ten computers where the entire work load is divided in an manual fashion. Each of the computers generate report based on their computational workflow which is synced to a database in a local server. The final report of the simulation is compiled using a python script from the individual work of each computer.}
	 \end{flushleft}
	 
	 This big computational work is manually divided in a fleet of ten computers each of which use implement three cores parallel processing to to complete their individual workload. All the computers used in the project are of the same spec which is illustrated in Table \ref{table:specs}. The report generated by each of the computer is synced to a database in a local server. The final report is generated automatically from the individual reports using a python script.\\ 
	    
    All the computation of this paper are done using Python 3.11.1. The technical specification of the computers used for the computation is presented in Table \ref{table:specs}.
    \begin{center}
            \captionof{table}{System Specification for the Experiment.}
            \begin{tabular}{|l|l|}
            	\hline
            	\label{table:specs}
                \textbf{Component} & \textbf{Spec}                                     \\ \hline
            	RAM & 16 GB                                             \\ \hline
            	CPU & i7-12700; 12   cores; 20 threads; 2.1 GHz  \\ \hline
            	Cache& 25 MB                                             \\ \hline
            	Graphics& Intel UHD Graphics 770  with shared Graphics Memory \\ \hline
            	Disk Space                                                                                      & 1 TB                                              \\ \hline
            \end{tabular}

    \end{center}

    \section{Result}\label{sec:result}
	Both the algorithms are tested against five different datasets of RSA numbers. In dataset 1, digits of length from 8 to 22 are adopted from \cite{rutkowski2020cryptanalysis} and the rest of the numbers are generated using the RSA number generator presented in \cite{rsagen}. For each number, we run 30 different instances of an algorithm and the reported average success is the mean (over the all the datasets) of the number of times the algorithm succeeded of the 30 instances tried out.
	
	\subsection{Simple Genetic Algorithm}\label{subsec:simplegaresult}
	Based on preliminary tests, we observe that the crossover rate of $ 50\% $ and mutation rate of $ 100\% $ works best for Simple GA. The detailed result of simple GA is illustrated in appendix \ref{app:simple_ga_output}. The biggest number it was able to factorize was 17 digit with an average success rate of $ 72\% $ as illustrated in Table \ref{tab:simple_ga_result}. This proposed GA is better that the GA proposed in \cite{rutkowski2020cryptanalysis} whose best best result was able to factorize a 17 digit number with an  average success rate of $ 3.33\% $. This shows that the space optimization scheme proposed in the paper has greatly improved the success rate of the algorithm.\\
	
	\begin{center}
		\begin{threeparttable}
			\caption{Average Result of Simple GA.}\label{tab:simple_ga_result}
			\begin{tabular}{lm{15mm} m{15mm} m{15mm} m{15mm} m{26mm} m{21mm}} 
				\rowcolor[HTML]{156082} 
				{\color[HTML]{FFFFFF} \textbf{Digits}} &
				{\color[HTML]{FFFFFF} \textbf{Success Rate(\%)}} &
				{\color[HTML]{FFFFFF} \textbf{Min Generation}} &
				{\color[HTML]{FFFFFF} \textbf{Max Generation}} &
				{\color[HTML]{FFFFFF} \textbf{Average Generation}} &
				{\color[HTML]{FFFFFF} \textbf{Time for 30 iterations(Seconds)}} &
				{\color[HTML]{FFFFFF} \textbf{Average time per iteration   (Seconds)}} \\
				\rowcolor[HTML]{C0E6F5} 
				8  & 100   & 0 & 1    & 1   & 26.80    & 0.89    \\
				9  & 100   & 0 & 1    & 1   & 32.88    & 1.10    \\
				\rowcolor[HTML]{C0E6F5} 
				10 & 100   & 0 & 4    & 1   & 52.16    & 1.74    \\
				11 & 100   & 0 & 34   & 7   & 407.27   & 13.58   \\
				\rowcolor[HTML]{C0E6F5} 
				12 & 100   & 0 & 36   & 6   & 276.85   & 9.23    \\
				13 & 100   & 0 & 228  & 35  & 1972.54  & 65.75   \\
				\rowcolor[HTML]{C0E6F5} 
				14 & 100   & 0 & 274  & 63  & 2362.38  & 78.75   \\
				15 & 100   & 0 & 1088 & 260 & 13461.85 & 448.73  \\
				\rowcolor[HTML]{C0E6F5} 
				16 & 92.67 & 0 & 1612 & 465 & 22268.45 & 742.28  \\
				17 & 72    & 0 & 1645 & 439 & 41021.92 & 1367.40
			\end{tabular}
			\begin{tablenotes}[para,flushleft]
				\footnotesize
				The table depicts average results of simple GA across all the five datasets over the which the simulations have been conducted. Detailed results over all the datasets can be found in the appendix \ref{app:simple_ga_output}.
			\end{tablenotes}
		\end{threeparttable}
	\end{center}
	
	In comparison with the molecular based algorithm proposed in \cite{mishra2016heuristic} and the firefly algorithm proposed in \cite{mishra2014multithreaded}, the simple GA proposed in this paper performs far better in terms of performance and accuracy. The biggest number the molecular based algorithm can factorize is a 14 digit (46 bit) number with a success rate of $ 69\% $ using on an average 2154.5 iterations. Similarly, in case of the same number the firefly algorithm achieved a success rate of $ 100\% $ with an average of 419 iterations. The simple GA proposed in this paper was able to factorize the same number with an success rate of $ 100\% $ with an average of 112 iterations. This is a significant development over the existing body of literature.
	
	\subsection{Sieve Algorithm}\label{subsec:sievealgoresult}
	Based on preliminary test, we observe that the crossover rate of $ 50\% $ and a mutation rate of $ 95\% $ works best for the sieve method. The detailed result of the sieve method is illustrated in the appendix \ref{app:sieve_algo_output}. It is clearly evident from the results of the simple GA and sieve method that the sieve method is often faster to converge to the solution in comparison to the simple GA. This has been possible due to the theoretical construct presented in equation \eqref{eq:sieve_method}. The algorithm has been able to factorize upto 23 digit (75 bit) number with and average success rate of $ 84\% $ and an average of $ 202 $ iterations as illustrated ins Table \ref{tab:sieve_algo_result}.\\
	
	\begin{center}
		\begin{threeparttable}
			\caption{Result of Sieve Method.}\label{tab:sieve_algo_result}
			\scriptsize
			\begin{tabular}{lm{15mm} m{15mm} m{15mm} m{15mm} m{26mm} m{21mm}}
				\rowcolor[HTML]{156082} 
				{\color[HTML]{FFFFFF} \textbf{Digits}} &
				{\color[HTML]{FFFFFF} \textbf{Success Rate(\%)}} &
				{\color[HTML]{FFFFFF} \textbf{Min Generation}} &
				{\color[HTML]{FFFFFF} \textbf{Max Generation}} &
				{\color[HTML]{FFFFFF} \textbf{Average Generation}} &
				{\color[HTML]{FFFFFF} \textbf{Time for 30 iterations (Seconds)}} &
				{\color[HTML]{FFFFFF} \textbf{Average time per iteration   (Seconds)}} \\
				\rowcolor[HTML]{C0E6F5} 
				8  & 100   & 0 & 1    & 0   & 21.77    & 0.73   \\
				9  & 100   & 0 & 1    & 0   & 20.53    & 0.68   \\
				\rowcolor[HTML]{C0E6F5} 
				10 & 100   & 0 & 1    & 0   & 21.12    & 0.70   \\
				11 & 100   & 0 & 8    & 1   & 77.92    & 2.60   \\
				\rowcolor[HTML]{C0E6F5} 
				12 & 100   & 0 & 10   & 1   & 75.50    & 2.52   \\
				13 & 100   & 0 & 48   & 6   & 483.77   & 16.13  \\
				\rowcolor[HTML]{C0E6F5} 
				14 & 100   & 0 & 123  & 14  & 795.38   & 26.51  \\
				15 & 100   & 0 & 224  & 38  & 2163.33  & 72.11  \\
				\rowcolor[HTML]{C0E6F5} 
				16 & 100   & 0 & 392  & 53  & 3368.85  & 112.30 \\
				17 & 98    & 0 & 1307 & 158 & 10081.22 & 336.04 \\
				\rowcolor[HTML]{C0E6F5} 
				18 & 93.33 & 0 & 1192 & 106 & 15097.72 & 503.26 \\
				19 & 94.67 & 0 & 685  & 60  & 5582.93  & 186.10 \\
				\rowcolor[HTML]{C0E6F5} 
				20 & 94.67 & 0 & 857  & 130 & 17159.60 & 571.99 \\
				21 & 96.67 & 0 & 1346 & 132 & 9496.42  & 316.55 \\
				\rowcolor[HTML]{C0E6F5} 
				22 & 83.33 & 0 & 1348 & 151 & 26497.59 & 883.25 \\
				23 & 84    & 0 & 1742 & 202 & 27797.82 & 926.59
			\end{tabular}
			\begin{tablenotes}[para,flushleft]
				\footnotesize
				The table depicts average results of sieve algorithm across all the five datasets over the which the simulations have been conducted. Detailed results over all the datasets can be found in the appendix \ref{app:sieve_algo_output}.
			\end{tablenotes}
		\end{threeparttable}
	\end{center}
	
	This is better in comparison with the existing literature where the best performing algorithm is ``Chromosome is m" which has been able to factorize is a number of 19 digit with a success rate of $ 3.33\% $.
	
	\subsection{Comparative Performance Analysis}
	
	The comparative performance analysis section will be discussed in two segments. Initially, we will compare the performance of the proposed algorithms with the best performing algorithm in the literature i.e. ``Chromosome is m". Finally, we will compare and contract between the performance of the proposed algorithms in this paper.
	
	\subsubsection{Performance Analysis with Literature} 
	In this section we will discuss the performance of Sieve algorithm and simple GA proposed in this paper with the best performing algorithm in the literature i.e. ``Chromosome is m" with respect to dataset 1 of the study, as the performance of the ``Chromosome is m" is available in the study \cite{rutkowski2020cryptanalysis}.\\
	
	In terms of success rate, both the proposed algorithm has outperformed ``Chromosome is m" as illustrated in figure \ref{fig:success_rate_comp}. It has been seen that, simple GA and sieve method has obtained $26.50\%$ and $321.89\%$ increment in success rate over the ``Chromosome is m" method respectively. It is evident from the figure \ref{fig:success_rate_comp} that, both the algorithms have higher success rate for larger number in comparison of the ``Chromosome is m" methods.

	\begin{flushleft}
		\includegraphics[width=\textwidth]{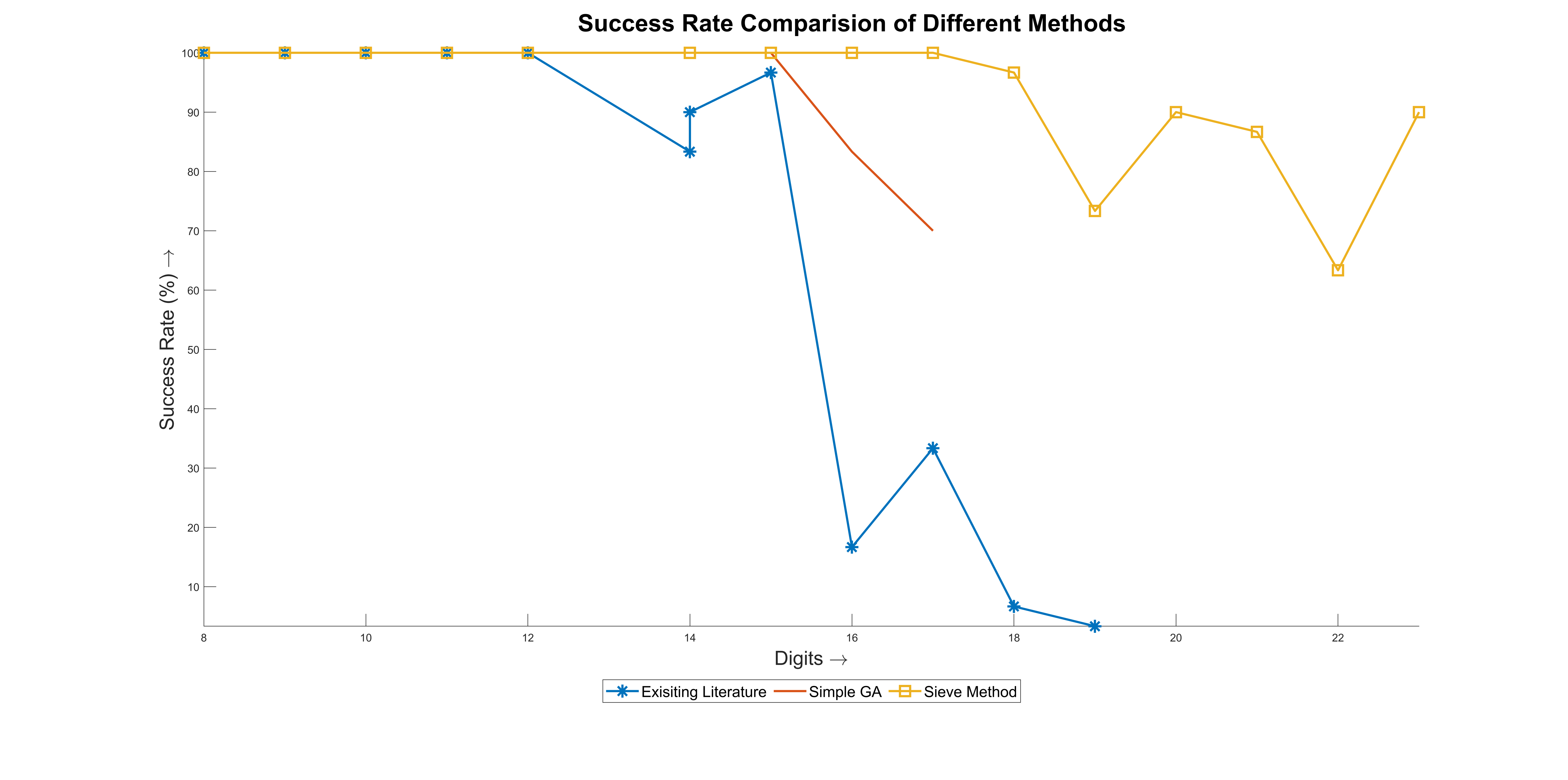}
		\captionof{figure}{Success rate comparison of  the proposed algorithms with the best algorithm in the existing literature.}\label{fig:success_rate_comp}
		\footnotesize{The figure depicts the success rate achieved by the three algorithms vs the digits of primes it has attempted to factorize. It is evident that the best performing algorithm is the Sieve method followed by the Simple GA. Both the algorithms proposed in this paper has out performed the best algorithm in the existing literature.}
	\end{flushleft}

	For ease of interpretation, a reduced number of maximum generations required for an algorithm to factorize a number indicates a more rapid convergence of the algorithm, thereby signifying superior performance. In terms of maximum iterations required for an algorithm to converge, sieve method outperforms all the other algorithms at almost all points as illustrated in figure \ref{fig:max_gen_rate_comp}. It has been observed that, simple GA and sieve method has obtained $41.91\%$ and $64.06\%$ decrement in maximum generation required for the algorithm to converge over the ``Chromosome is m" method respectively. It is evident from the figure \ref{fig:max_gen_rate_comp} that, both the algorithms converges faster for larger number in comparison of the ``Chromosome is m" methods given the lower max generation required.

		\begin{flushleft}
		\includegraphics[width=\textwidth]{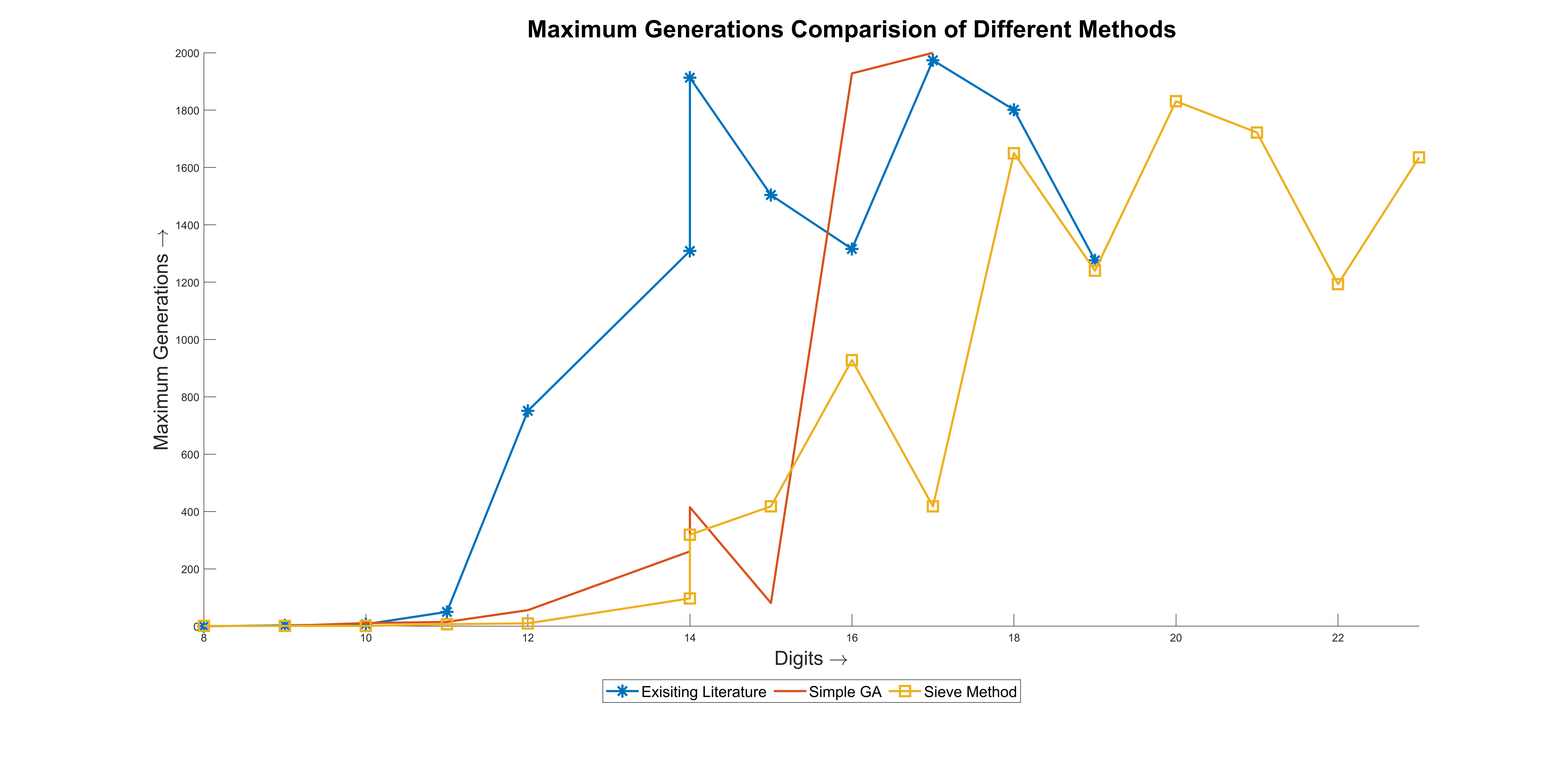}
		\captionof{figure}{Max generation comparison of  the proposed algorithms with the best algorithm in the existing literature.}\label{fig:max_gen_rate_comp}
		\footnotesize{The figure depicts the maximum generation required for the three algorithms to converge vs the digits of primes it has attempted to factorize. For ease of interpretation, the lower the maximum generation required, the faster the algorithm succeed to factorize the number, i.e. the flatter the curve the better. It is evident that the best performing algorithm is the Sieve method followed by the Simple GA.}
	\end{flushleft}
	
	\subsubsection{Performance Analysis of the Proposed Model}
	
	This section discusses the performance comparison of the two proposed algorithms in this study namely the simple GA and sieve method. The performance comparison is done based on the simulation of the five datasets. Summarized result for each of the algorithms has been presented in Table \ref{tab:simple_ga_result} and \ref{tab:sieve_algo_result} and detailed result of the simulations has been presented in appendix \ref{app:simple_ga_output} and \ref{app:sieve_algo_output} respectively. We shall analyze the performance based on four KPIs namely: success rate, maximum generation, average generation and average time per iteration. The comparison of the two algorithms based on these four KPIs are well illustrated in figure \ref{fig:performance_comp}.\\

	\begin{flushleft}
		\includegraphics[width=\textwidth]{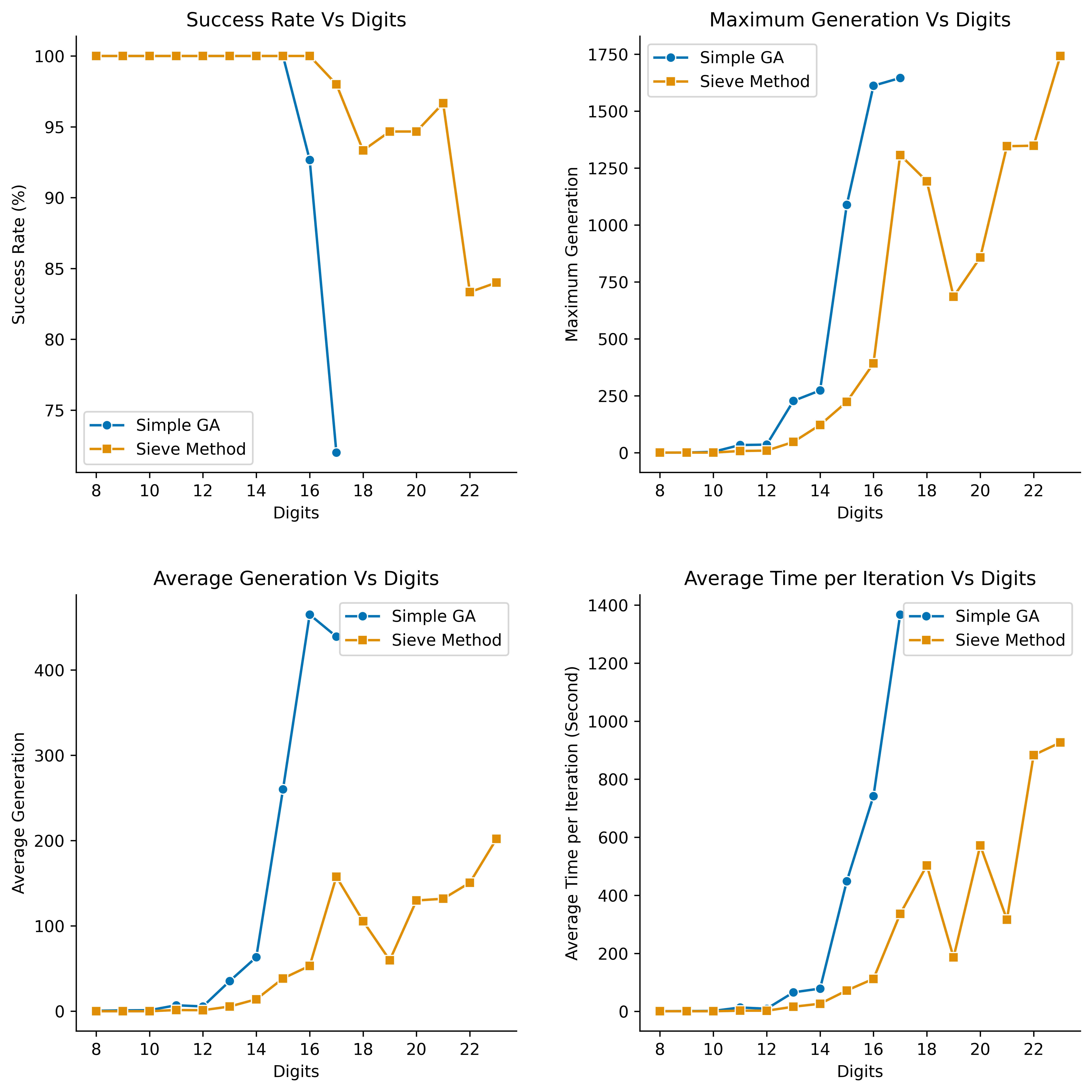}
		\captionof{figure}{Performance comparison of the proposed algorithms.}\label{fig:performance_comp}
		\footnotesize{The figure depicts the performance comparison of the proposed model namely simple GA and sieve method in terms of success rate, maximum generation, average generation, and average time per iteration. It is evident from the figure that the sieve method is an optimized algorithm with respect to the simple GA in terms of the four KPIs.}
	\end{flushleft}
	
	In terms of success rate, the it can be seen from figure \ref{fig:performance_comp} that the sieve method has higher success rate for larger number in comparison to simple GA method and for smaller numbers the algorithm is on par with simple GA if not better.\\
	
	In terms of maximum generations required for the algorithm to converge, sieve method performs better with respect to the simple GA algorithm for all possible choices of digit length. Sieve methods has been found to provide around $55.33\%$ decrement of maximum generation on average in comparison to simple GA.\\
	
	In terms of average generation required for the algorithms to converge, the sieve method provides better result in comparison to the simple GA method given it's flatter curve as the number of digits increases in comparison to the simple GA. Simulation results shows that sieve method provides around $85.87\%$ decrement in average generation on average in compression to simple GA.\\
	
	Finally in terms of average time per iteration, sieve methods has exhibits $65.55\%$ decrement on average with respect to simple GA method which implies faster convergence of the sieve method.\\

	From the aforesaid discussion, we can easily conclude that the sieve method is significantly better than the simple GA and both the algorithms proposed in this paper provide significantly better result in comparison to the best algorithm in the existing literature i.e. ``Chromosome is m".

    \section{Future Scope of the Study}\label{sec:lim}
    
    The results of the study are very promising  and has shown a great deal of improvement over the existing literature. There some avenues that can be explored to push the outcomes of this study beyond it's current capabilities.\\
    
    The best algorithm proposed in this paper i.e. the sieve method relies heavily on two instances of GA one working top of another to optimize the selection of a,d of equation \eqref{eq:sieve_method} from a discrete range of values to effectively shrink the search space for faster convergence. This slows down the execution speed as the number of digits increases hence require higher computational power and time. This is problematic as GA is very sensitive to hyper parameter tuning. Relation between the digit length of the big semi-prime, $M$ with that of the sieve (a, d) of it's prime factor is an avenue that need further research. If a theoretical construct can be developed bridging the afore mentioned relation, them the sieve method can be reconstructed based on a single instance of GA rather than two which can make the algorithm much more potent to factorize bigger semi-primes with very little computational power compared that of it's present construct.\\
    
    The main reason behind the results achieved in this paper over the existing literature is because of the  search space optimization schemes adopted based on the theoretical construct of the semi-primes and primes. By construct, the algorithms proposed in this paper targets the smallest prime factor which is not necessarily at the central region of the search space. Research on techniques to effectively narrow down the location of the prime factor in the search space based on the construct of the big semi-prime $M$ can effectively boost up the speed and accuracy of the algorithm.

    \section{Conclusion }\label{sec:conclusion}
	Application of computational intelligence approach in the filed of cryptanalysis is a fairly new avenue of research which will gain more momentum in the near future with the rise of computational power. This paper initially proposes a conjecture on the distribution of digits in primes as length of prime tends to infinity. The algorithms presented in this paper namely simple GA and sieve method are improvement over the existing literature which exploits the search space shrinkage schemes and theoretical constructs of prime numbers to generate  $26.50\%$, $321.89\%$ increment in success rate and $41.91\%$ and $64.06\%$ decrement in maximum generation required for the algorithm to converge over the existing literature respectively which is a significant improvement. The largest number the algorithms proposed in this paper can factor are 23 digits (75 bits) long with an average success rate of $84\%$. An important note in this regard is that although the paper showed promising result over the existing literature, the biggest number the algorithm can factorize using a general user grade computational environment is a lot smaller than the current industry standard of RSA cryptosystem hence the proposed development is not an eminent threat to industry standard security system anywhere. The development in this avenue of research is fairly recent but with the rise of computational power in the near future, this class of algorithms can be potent threat to cryptosystem.
	
	\section*{Acknowledgments}
	The research by M. Kamrujjaman was partially supported by the University Grants Commission (UGC), 
	and  the  University of Dhaka, Bangladesh.
	\section*{Conflict of interest}
	The authors declare no conflict of interest.
	
	\section*{Declaration of generative AI and AI-assisted technologies in the writing process}
	During the preparation of this work the author(s) never used any tools/service of AI and AI-assisted technologies.
	
	\section*{Ethical approval}
	No consent is required to publish this manuscript.
	
	\section*{Author contributions}
	Conceptualization,  MAM and MK; methodology, MAM; software, MAM; validation, MK; formal analysis, MAM; investigation, MK; resources, MK and MAM;
	original draft preparation, MAM and MK; review and editing,  MK; supervision, MK. All authors have read and agreed to the published version of the manuscript.

    \newpage
    \appendix
    \raggedright \textbf{\Huge Appendix }

    \section{Results of Sieve Algorithm Simulation}\label{app:sieve_algo_output}
	\begin{center}
		\includegraphics[page=1,height=0.85\textheight,width=\textwidth]{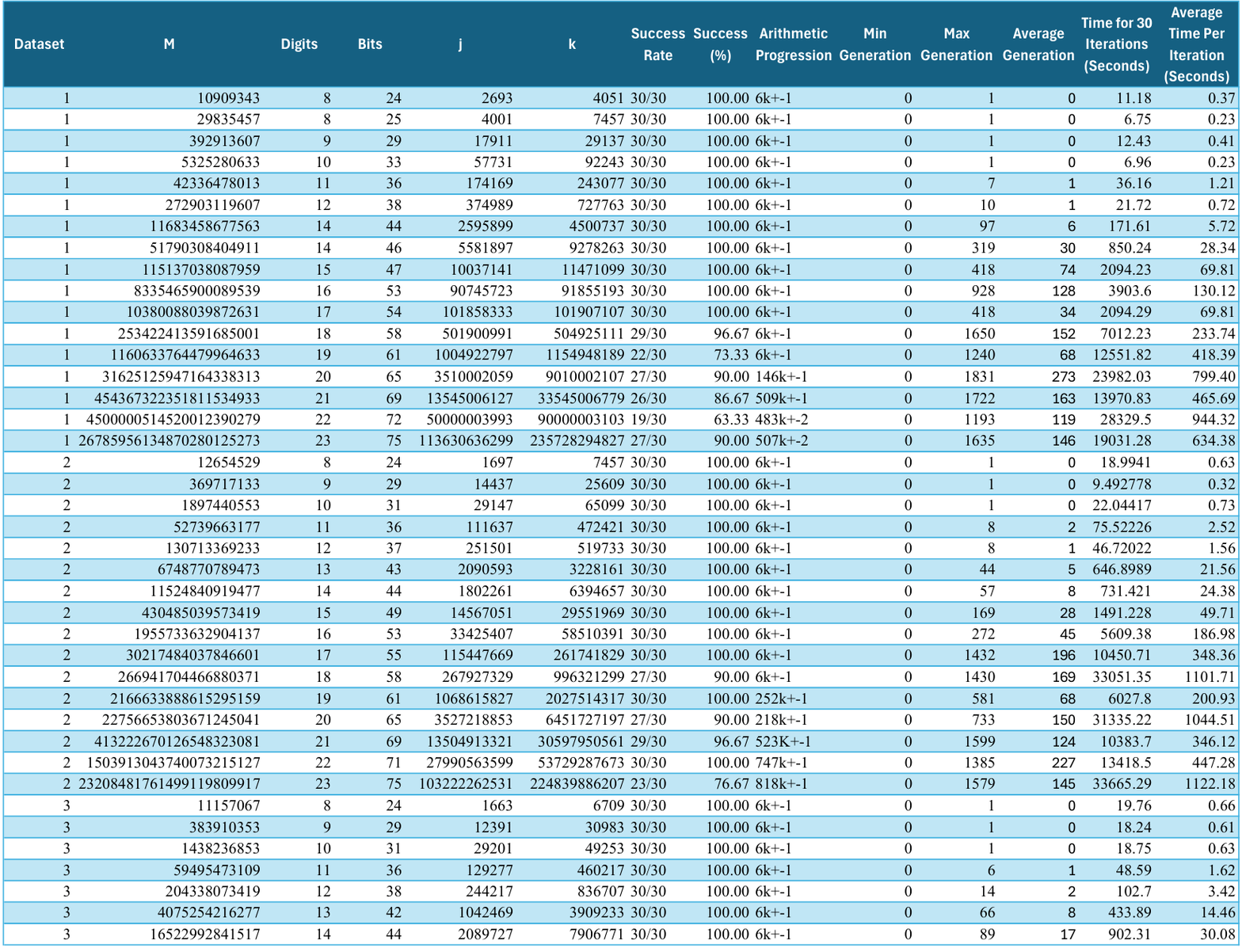}
		\includegraphics[page=2,height=\textheight,width=\textwidth]{sieve_output.pdf}
		\includegraphics[page=3,height=0.25\textheight,width=\textwidth,trim={0 6inch 0 0},clip]{sieve_output.pdf}
	\end{center}
	
	\section{Results of Simple GA Simulation}\label{app:simple_ga_output}
	\begin{center}
		\includegraphics[page=1,height=0.6\textheight,width=\textwidth]{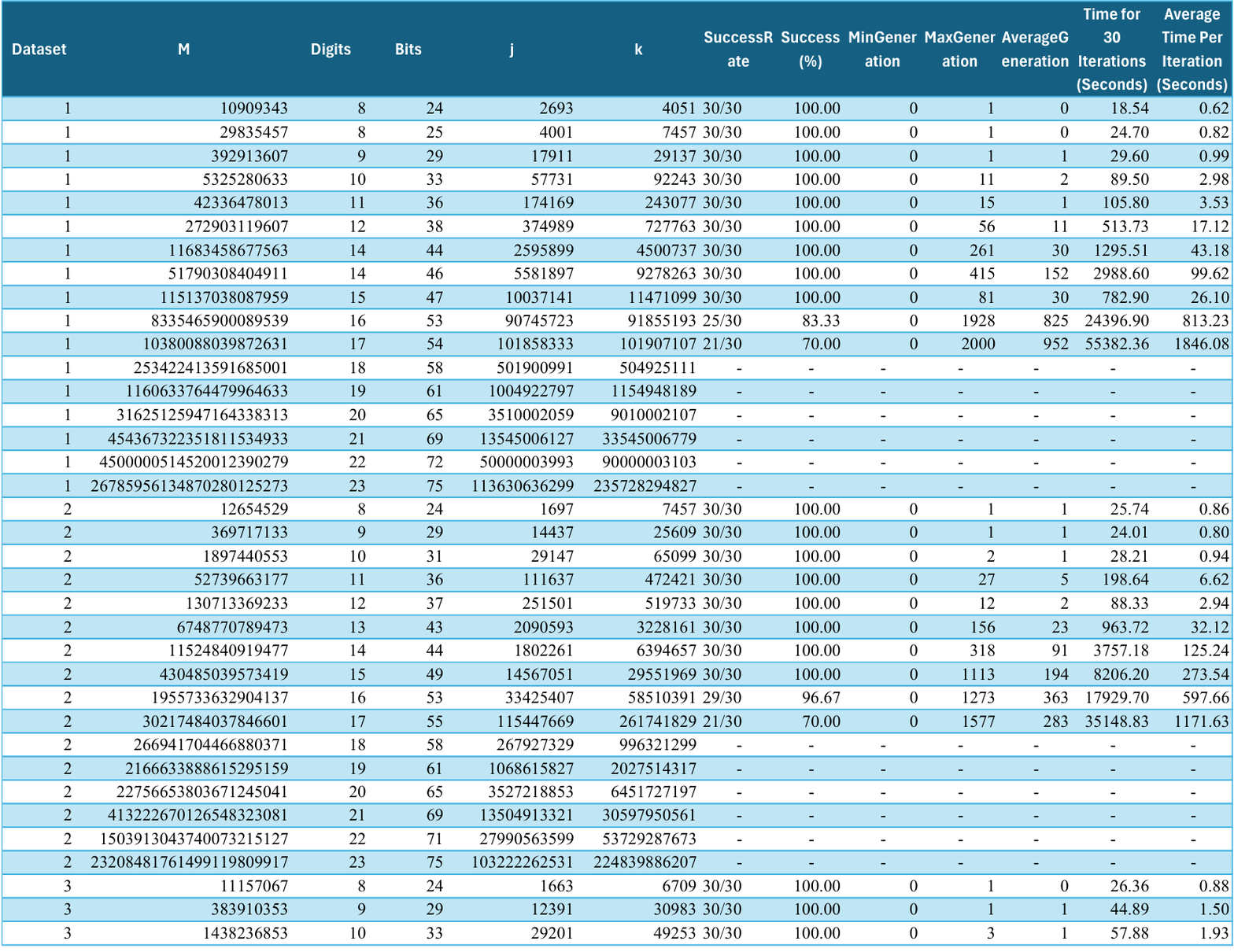}
		\includegraphics[page=2,height=\textheight,width=\textwidth]{simple_ga_output.pdf}
		\includegraphics[page=3,height=0.25\textheight,width=\textwidth,trim={0 5inch 0 0},clip]{simple_ga_output.pdf}
	\end{center}

    \nocite{*}
    \bibliographystyle{unsrt}
    \bibliography{reference}

\begin{thebibliography}{10}

\bibitem{rivest1978method}
Ronald~L Rivest, Adi Shamir, and Leonard Adleman.
\newblock A method for obtaining digital signatures and public-key
  cryptosystems.
\newblock {\em Communications of the ACM}, 21(2):120--126, 1978.

\bibitem{mishra2016heuristic}
Mohit Mishra, Utkarsh Chaturvedi, and Kaushal~K Shukla.
\newblock Heuristic algorithm based on molecules optimizing their geometry in a
  crystal to solve the problem of integer factorization.
\newblock {\em Soft Computing}, 20(9):3363--3371, 2016.

\bibitem{yampolskiy2010application}
Roman~V Yampolskiy.
\newblock Application of bio-inspired algorithm to the problem of integer
  factorisation.
\newblock {\em International Journal of Bio-Inspired Computation},
  2(2):115--123, 2010.

\bibitem{mudgal2017application}
Piyush~Kumar Mudgal, Rajesh Purohit, Rajesh Sharma, and Mahendra~Kumar Jangir.
\newblock Application of genetic algorithm in cryptanalysis of mono-alphabetic
  substitution cipher.
\newblock In {\em 2017 International Conference on Computing, Communication and
  Automation (ICCCA)}, pages 400--405. IEEE, 2017.

\bibitem{shikhare2015cryptanalysis}
Aparna Shikhare.
\newblock Cryptanalysis of the purple cipher using random restarts.
\newblock 2015.

\bibitem{brown2009genetic}
Joseph~Alexander Brown, Sheridan Houghten, and Beatrice Ombuki-Berman.
\newblock Genetic algorithm cryptanalysis of a substitution permutation
  network.
\newblock In {\em 2009 IEEE Symposium on Computational Intelligence in Cyber
  Security}, pages 115--121. IEEE, 2009.

\bibitem{ferriman2014solving}
Benjamin Ferriman and Charlie Obimbo.
\newblock Solving for the rc4 stream cipher state register using a genetic
  algorithm.
\newblock {\em International Journal of Advanced Computer Science and
  Applications}, 5(5), 2014.

\bibitem{ma2011evolutionary}
Eddie Yee-Tak Ma and Charlie Obimbo.
\newblock An evolutionary computation attack on one-round tea.
\newblock {\em Procedia Computer Science}, 6:171--176, 2011.

\bibitem{ribaric2017genetic}
Tim Ribaric and Sheridan Houghten.
\newblock Genetic programming for improved cryptanalysis of elliptic curve
  cryptosystems.
\newblock In {\em 2017 IEEE Congress on Evolutionary Computation (CEC)}, pages
  419--426. IEEE, 2017.

\bibitem{chan2002automatic}
David~Michael Chan.
\newblock Automatic generation of prime factorization algorithms using genetic
  programming.
\newblock {\em Genetic Algorithms and Genetic Programming at Stanford}, pages
  52--57, 2002.

\bibitem{chang2005fast}
Weng-Long Chang, Minyi Guo, and MS-H Ho.
\newblock Fast parallel molecular algorithms for dna-based computation:
  factoring integers.
\newblock {\em IEEE Transactions on Nanobioscience}, 4(2):149--163, 2005.

\bibitem{meletiou2002first}
G~Meletiou, DK~Tasoulis, Michael~N Vrahatis, et~al.
\newblock A first study of the neural network approach to the rsa cryptosystem.
\newblock In {\em IASTED 2002 Conference on Artificial Intelligence}, pages
  483--488. IASTED Calgary, 2002.

\bibitem{mishra2014multithreaded}
Mohit Mishra, Utkarsh Chaturvedi, and Saibal~K Pal.
\newblock A multithreaded bound varying chaotic firefly algorithm for prime
  factorization.
\newblock In {\em 2014 IEEE International Advance Computing Conference (IACC)},
  pages 1322--1325. IEEE, 2014.

\bibitem{rutkowski2020cryptanalysis}
Emilia Rutkowski and Sheridan Houghten.
\newblock Cryptanalysis of rsa: Integer prime factorization using genetic
  algorithms.
\newblock In {\em 2020 IEEE Congress on Evolutionary Computation (CEC)}, pages
  1--8. IEEE, 2020.

\bibitem{xingbo2019number}
WANG Xingbo.
\newblock Number of digits in two integers and their multiplication.
\newblock {\em Journal of Advances in Applied Mathematics}, 4(2), 2019.

\bibitem{rsagen}
Big primes \url{https://bigprimes.org/RSA-challenge}.
\newblock Accessed on 05 April 2024.

\end{thebibliography}
\end{document}